\documentclass[12pt]{amsart}
\usepackage{amsmath,amssymb,amsfonts, amscd}
\usepackage{enumerate}
\usepackage{anysize}
\marginsize{3cm}{3cm}{3cm}{3cm}
\input xy
\xyoption{all}
\usepackage{pb-diagram}
\usepackage[all]{xy}
\input xy
\xyoption{all}

\theoremstyle{plain}
\newtheorem{thm}{Theorem}
\newtheorem{theorem}{Theorem}[section]

\newtheorem{lemma}{Lemma}
\newtheorem{cor}[theorem]{Corollary}

\newtheorem{conj}{Conjecture}

\theoremstyle{definition} \theoremstyle{definition}

\theoremstyle{remark}

\newcommand{\GG}{\mathbb{G}}

\newcommand{\Q}{\mathbb{Q}}

\newcommand{\U}{\mathcal{U}}
\newcommand{\Z}{\mathbb{Z}}

\newcommand{\R}{\mathbb{R}}

\newcommand{\C}{\mathbb{C}}

\newcommand{\wT}{\widehat{T}}

\newcommand{\wrho}{\widehat{\rho}}

\def\SL{{\rm SL}}

\def\Sp{{\rm Sp}}
\def\St{{\rm St}}

\def\U{{\rm U}}

\def\GL{{\rm GL}}
\def\PGL{{\rm PGL}}

\def\SO{{\rm SO}}

\begin{document}

\title[A  character relationship on $\GL_n$]
{A  character relationship on $\GL_n$}

\author{Dipendra Prasad}

\address{School of Mathematics, Tata Institute of Fundamental
Research, Colaba, Mumbai-400005, INDIA}
\email{dprasad@math.tifr.res.in}
\maketitle

\begin{abstract}In this paper we consider the character of a finite dimensional algebraic 
representation of $\GL_{mn}(\C)$ restricted to a particular disconnected component 
of the normalizer of the Levi subgroup $\GL_m(\C)^n$ of $\GL_{mn}(\C)$, generalizing
a theorem of Kostant on the character values at the Coxeter element. 
\end{abstract}

\section{Introduction}

Studying representation theory of disconnected groups in the context of real and 
$p$-adic groups has been an important topic of study, finding  impressive 
applications such as to the theory of basechange, cf. [La],  which has been a key instrument in 
all recent proofs of reciprocity theorems in number theory, and eventually to the Fermat's Last theorem!

Representation theory of disconnected 
algebraic groups can also be studied in a similar vein, see for instance the paper [KLP].
In this paper, we prove a simple result on the restriction of representations of  
connected reductive algebraic group (actually we are able to handle  only $\GL_n(\C)$ here) 
to suitable disconnected algebraic groups whose
connected component of identity is a Levi subgroup with `large enough' normalizer.

We begin with a very beautiful theorem of Kostant, Theorem 2 in [Ko],  which is about the character of a finite dimensional representation of a  semi-simple algebraic group $G$ at 
the Coxeter conjugacy class. Recall that one usually defines a Coxeter element --- or, rather a conjugacy class --- 
in a Weyl group (as a product of simple reflections), in this case in $N(T)/T$, where $T$ is a maximal torus in $G$, with $N(T)$ its normalizer in $G$.
An arbitrary lift of this conjugacy class in $N(T)/T$ to $N(T)$ gives a well-defined conjugacy class in $G$ which  we will denote by $c(G)$. 
Here is the theorem of Kostant.

\begin{thm} Let $G$ be a semi-simple algebraic group over $\C$, and $\pi$ a finite dimensional irreducible representation of $G$. 
Then the character $\Theta_\pi$ of $\pi$ at the element $c(G)$ takes one of the values $1,0,-1$.
\end{thm}

\section{The main theorem}

Before we  turn to the main theorem of this paper,  we note   the following well-known and 
elementary lemma whose proof we will omit. 
(This lemma replaces the  observation  that an arbitrary lift to $N(T)$ of a Coxeter   element in   $N(T)/T$ belongs to a unique
conjugacy class in $G$.)

\vspace{4mm}

\begin{lemma}Let $G$ be an arbitrary group. Consider the semi-direct product 
$G^n \rtimes \Z/n$ 
in which $\Z/n$ operates on $G^n$ by the $n$-cycle: $(g_1,\cdots, g_{n-1},g_n)
\rightarrow (g_n,g_1,\cdots,g_{n-1})$; we call this generator of $\Z/n$ by $\sigma$. Then 
there exists a bijective correspondence between conjugacy class of elements of
$G^n \rtimes \Z/n$  of the form $(g_1,\cdots, g_{n-1},g_n) \rtimes \sigma$ 
and conjugacy class of elements of $G$ given by the Norm mapping, ${\rm Nm}: G^n \rtimes \Z/n \rightarrow  G$, given by: 
$$(g_1,\cdots, g_{n-1},g_n) \rtimes \sigma \in G^n \rtimes \Z/n \rightarrow g_1\cdots g_n \in G.$$
\end{lemma}

 \vspace{4mm}

The aim of this paper is to calculate the character of irreducible representations of 
$\GL_{mn}(\C)$ on the subgroup ${\GL_m(\C)}^n \rtimes \Z/n$ 
which sits
naturally inside $\GL_{mn}(\C)$ at  elements of the subgroup which have projection a  fixed
generator $\sigma$ of $\Z/n$ generalizing the theorem of Kostant recalled above for $\GL_n(\C)$. Since elements of $\GL_m(\C)$ have $n$-th roots, it suffices by the previous lemma to calculate the character at very special elements, those of the form
$(g,\cdots, g,g) \rtimes \sigma .$ This is what we shall do here.

Consider the natural embedding of groups,
$$\iota: \GL_m(\C) \times \GL_n(\C) \longrightarrow \GL_{mn}(\C),$$ 
obtained by taking the tensor product of vector spaces over $\C$ of dimensions $m,n$.
Thus, we need to calculate the characters of irreducible representations of $\GL_{mn}(\C)$ at the elements of the form $\iota(t\cdot c_n)$, which we will simply write as $t \cdot c_n$,
 where $t$ is any 
element of $\GL_m(\C)$, and $c_n = c(\GL_n(\C))$.

\begin{thm} Let $\pi$ be an irreducible finite dimensional representation of $\GL_{mn}(\C)$ with character $\Theta_\pi$ and with highest weight 
$$\underline{\lambda}= \lambda_1 \geq \lambda_2 \geq \cdots \geq \lambda_{mn-1}\geq \lambda_{mn}.$$
Let $\underline{\delta}_{mn}$ be defined to be the $mn$-tuple of integers,
$$mn-1 \geq mn-2 \geq \cdots \geq 1 \geq 0.$$
Then the character  $\Theta(t \cdot c_n)$ is not identically zero if and only 
if the $mn$-tuple of integers appearing in $\underline{\lambda} + \underline{\delta}_{mn}$ 
represent each residue class in $\Z/n$ exactly $m$-times.
Assume this to be the case. Then for each residue class in $\Z/n$ represented by 
$0\leq i < n$, let $\pi_i$ be the representation of $\GL_m(\C)$ with highest weight 
$\underline{\mu}_i$ such that $\underline{\mu}_i+\underline{\delta}_m$ are 
the integers (exactly $m$ of them by hypothesis) among, 
$$[\underline{\lambda} + \underline{\delta}_{mn} -i]/n.$$ 

 Then for the  irreducible representations
$\pi_1,\cdots, \pi_n$,  of $\GL_m(\C)$ with highest weights $\underline{\mu}_1,\cdots, 
\underline{\mu}_n$, and 
with characters $\Theta_1,\cdots,\Theta_n$,
$$\Theta(t \cdot c_n)
= \pm \Theta_1(t^n) \Theta_2(t^n)\cdots \Theta_n(t^n).$$ 
\end{thm}

\vspace{4mm}

{\bf Remark :} Before we begin the proof of the theorem, it may be worth pointing out that 
the representations $\pi_i$ of $\GL_m(\C)$ occurring in this theorem through the character 
identity 
$$\Theta(t \cdot c_n)
= \pm \Theta_1(t^n) \Theta_2(t^n)\cdots \Theta_n(t^n).$$ 
are well defined only 
up to  twisting by characters of $\C^\times$ via the determinant map from $\GL_m(\C)$ to 
$\C^\times$ (with the only constraint that the product of these twisting characters 
is trivial). 
This character identity implies that if $z_\pi$ is the 
central character of $\pi$ (an integer by which the center of $\GL_{mn}(\C)$ which is 
$\C^\times$  operates 
on $\pi)$, and if $z_i$ are the central characters of $\pi_i$, then we must have,
$$z_\pi = z_1^n\cdot z_2^n\cdots z_n^n.$$
Since the central character of $\pi$ is given by the integer $\sum_i \lambda_i$, and similarly of $\pi_i$, that the above identity of central characters  holds in our case follows from the following identity:
$$ \sum_{i=0}^{i=mn-1} i = m\sum_{i=0}^{i=n-1} i + n^2\sum_{i=0}^{i=m-1} i.$$

\vspace{4mm}

\begin{proof}(of the Theorem) We begin by proving that if the $mn$-tuple of integers represented by 
$\underline{\lambda} + \underline{\delta}_{mn}$ do not represent
each residue class in $\Z/n$ by exactly $m$-integers, then the character value is
identically zero at the elements of the form $t\cdot c_n$, where we take 
$t$ to be the diagonal matrix in $\GL_m(\C)$ with diagonal entries $(t_1,t_2,\cdots, t_n)$, and let $c_n$ be the diagonal matrix 
in $\GL_n(\C)$ consisting of the $n$-th roots of unity all with multiplicity one.

Thus $t\cdot c_n$ is the diagonal matrix in $\GL_{mn}(\C)$ 
with entries $t_i\omega^j$ where 
$\omega$ is a primitive $n$-th root of unity in $\C$, and 
$1 \leq i \leq m$, and $0\leq j < n$.

Under the assumption on  the $mn$-tuple of integers represented by 
$\underline{a}_i= \underline{\lambda} + \underline{\delta}_{mn}$, by re-ordering the 
indices, assume that the integers $a_1,a_2,\cdots, a_{m+1}$ have the same residue modulo $n$.

Now we write out a part of the matrix whose determinant represents the numerator of the
Weyl character formula:

$$
\left ( \begin{array}{lllllllllll} 
t_1^{a_1} &  t_2^{a_1}   & \cdots  & t_m^{a_1} & 
\omega_1t_1^{a_1}  &  \cdots &   \omega_1t_m^{a_1} &    
\omega_2t_1^{a_1}   &  \cdots &   \omega_2t_m^{a_1} &    
\cdots \\
t_1^{a_2} &  t_2^{a_2}  &  \cdots & t_m^{a_2} & 
\omega_1t_1^{a_2}    & \cdots &   \omega_1t_m^{a_2} & 
\omega_2t_1^{a_2}    & \cdots &   \omega_2t_m^{a_2} & 
\cdots \\
\cdot  &  \cdot   &   \cdots   &  \cdot  
&  \cdot  &  \cdots  & \cdot   &  \cdot  &  \cdots  & \cdot  & \cdot  \\
\cdot  &  \cdot   &   \cdots   &  \cdot  
&  \cdot  &  \cdots  & \cdot   &  \cdot  &  \cdots  & \cdot  & \cdot  \\
t_1^{a_{m+1}} &  t_2^{a_{m+1}}   & \cdots &  t_m^{a_{m+1}} & 
\omega_1t_1^{a_{m+1}}  & \cdots &  \omega_1t_m^{a_{m+1}} &
\omega_2t_1^{a_{m+1}}    & \cdots &  \omega_2t_m^{a_{m+1}} 
& \cdots \\

\end{array}\right)\cdot$$

\vspace{3mm}

In this $(m+1) \times mn$-matrix, 
the first $(m+1) \times m$-matrix 
repeats itself after being scaled by $\omega_1$, 
an $n$-th root of unity, and then again  after being scaled by $\omega_2$, 
an $n$-th root of unity, and so on. Clearly the rank of this  $(m+1) \times mn$-matrix is 
$\leq m$, which proves that the determinant of the $mn \times mn$ matrix representing 
the numerator of the Weyl character is zero. The element $t\cdot c_n$  in $\GL_{mn}(\C)$ 
is regular for generic choice of $t$ in $\GL_m(\C)$, so the Weyl denominator is non-zero.
Thus the character of the representation $\pi$ at the element $t\cdot c_n$  in $\GL_{mn}(\C)$ 
is identically zero.

The proof of the character relationship in the theorem is by a direct manipulation with the character 
formula for an irreducible highest weight module of $\GL_{mn}(\C)$ given in terms of quotients of Vandermonde determinants, and depends on the observation that 

$$
\det \left ( \begin{array}{llll} 

X_1 & \omega X_1 & \cdots & \omega^{n-1} X_1  \\
 X_2 & \omega^2 X_2 & \cdots & \omega^{2(n-1)} X_2  \\

\cdot  & \cdot & \cdots & \cdot  \\

X_n & \omega^n X_n & \cdots & \omega^{n(n-1)} X_n  \\

\end{array}\right) = c \det(X_1) \det (X_2) \cdots \det (X_n),$$
where $\omega$ is a primitive $n$-th root of unity, 
$X_i$ are $m \times m$ matrices, and  $c=\prod_{i<j}(\omega^i-\omega^j)^m$.

We will not give a more detailed proof of the theorem in general,  
but show the details for  $m=n=2$.
Since the Coxeter element $c_2$ for $\GL_2(\C)$ is the diagonal matrix with entries $(1,-1)$, the theorem proposes to calculate the 
character of a representation of $\GL_4(\C)$ say with highest weight 
$\lambda_1 \geq \lambda_2 \geq \lambda_3\geq \lambda_4$ 
at the diagonal elements
of the form $(t_1,t_2,-t_1,-t_2)$. 

The proof of the theorem will be achieved by the explicit character of the representation $\pi$ given as the quotient of two 
Vandermonde determinants, the numerator of which is the determinant of the $4 \times 4$ matrix:

$$
\left ( \begin{array}{llll} 
t_1^{\lambda_1+3} &  t_2^{\lambda_1+3}  & (-t_1)^{\lambda_1+3} & (-t_2)^{\lambda_1+3} \\
t_1^{\lambda_2+2} &  t_2^{\lambda_2+2}  & (-t_1)^{\lambda_2+2} & (-t_2)^{\lambda_2+2} \\
t_1^{\lambda_3+1} &  t_2^{\lambda_3+1}  & (-t_1)^{\lambda_3+1} & (-t_2)^{\lambda_3+1} \\
t_1^{\lambda_4} &  t_2^{\lambda_4}  & (-t_1)^{\lambda_4} & (-t_2)^{\lambda_4} 
\end{array}\right)\cdot$$

It is easy to see that
for the determinant to be nonzero, it is necessary that out of 
$\lambda_1+3,  \lambda_2+2, \lambda_3+1, \lambda_4 $ two are even and two are odd. Assume 
for definiteness that 
$\lambda_1+3,  \lambda_2+2 $ are even, and
$ \lambda_3+1, \lambda_4 $ are odd, in which 
case, the Weyl numerator becomes,

$$
\left ( \begin{array}{llll} t_1^{\lambda_1+3} &  t_2^{\lambda_1+3}  & t_1^{\lambda_1+3} & 
t_2^{\lambda_1+3} \\
t_1^{\lambda_2+2} &  t_2^{\lambda_2+2}  & t_1^{\lambda_2+2} & t_2^{\lambda_2+2} \\
t_1^{\lambda_3+1} &  t_2^{\lambda_3+1}  & -t_1^{\lambda_3+1} & -t_2^{\lambda_3+1} \\
t_1^{\lambda_4} &  t_2^{\lambda_4}  & -t_1^{\lambda_4} & -t_2^{\lambda_4} \end{array}\right)\cdot$$

Since

$$\det \left ( \begin{array}{cccc} a_1 & a_2 & a_1 & a_2 \\
b_1 & b_2 & b_1 & b_2 \\
c_1 & c_2 & -c_1 & -c_2 \\
d_1 & d_2 & -d_1 & -d_2 \end{array} \right ) = 
4 \det \left ( \begin{array}{cccc} a_1 & a_2 & 0 & 0 \\
b_1 & b_2 & 0 & 0 \\
0 & 0 & -c_1 & -c_2 \\
0 & 0  & -d_1 & -d_2 \end{array} \right ),
$$  
we find that the Weyl numerator is the same as

\begin{eqnarray*}
&{}& 4 \det \left ( \begin{array}{ll}t_1^{\lambda_1+3} & t_2^{\lambda_1+3} \\
t_1^{\lambda_2+2} & t_2^{\lambda_2+2}\end{array}\right) \cdot \det \left( \begin{array}{ll}t_1^{\lambda_3+1} & t_2^{\lambda_3+1} \\
t_1^{\lambda_4} & t_2^{\lambda_4} \end{array}\right)  \\
& = &4 t_1t_2\det \left ( \begin{array}{ll}t_1^{\lambda_1+1} & t_2^{\lambda_1+1} \\
t_1^{\lambda_2} & t_2^{\lambda_2}\end{array}\right) \det \left ( \begin{array}{ll}t_1^{\lambda_3+2} & t_2^{\lambda_3+2} \\
t_1^{\lambda_4+1} & t_2^{\lambda_4+1}\end{array}\right) \\
& = &4t_1t_2\det \left ( \begin{array}{ll}(t_1^2)^{\frac{\lambda_1+1}{2}} & (t_2^2)^{\frac{\lambda_1+1}{2}} \\
(t_1^2)^{\frac{\lambda_2}{2}} & (t_2^2)^{\frac{\lambda_2}{2}}\end{array}\right)  \det \left ( \begin{array}{ll}(t_1^2)^{\frac{\lambda_3+2}{2}} & (t_2^2)^{\frac{\lambda_3+2}{2}} \\
(t_1^2)^{\frac{\lambda_4+1}{2}} & (t_2^2)^{\frac{\lambda_4+1}{2}}\end{array}\right).
\end{eqnarray*}

Further, the Weyl denominator at $(t_1,t_2,-t_1,-t_2)$ can be checked to be:
$$4t_1t_2(t_1^2-t_2^2)^2.$$

Therefore, 
\begin{eqnarray*}
\Theta(t)
& = & 
\frac{(t_1t_2)\cdot \det \left (\begin{array}{ll}(t_1^2)^{\frac{\lambda_1+1}{2}} & (t_2^2)^{\frac{\lambda_1+1}{2}} \\
(t_1^2)^{\frac{\lambda_2}{2}} & (t_2^2)^{\frac{\lambda_2}{2}}\end{array}\right) \cdot  \det \left (\begin{array}{ll}(t_1^2)^{\frac{\lambda_3+2}{2}} & (t_2^2)^{\frac{\lambda_3+2}{2}} \\
(t_1^2)^{\frac{\lambda_4+1}{2}} & (t_2^2)^{\frac{\lambda_4+1}{2}}\end{array}\right)}{t_1t_2(t_1^2-t_2^2)^2} \\
& = & \Theta_1(t^2)\cdot \Theta_2(t^2).
\end{eqnarray*}

\end{proof}

{\bf Remark 1:} 
The above theorem can be re-written for general elements of 
${\GL_m(\C)}^n \rtimes \Z/n$ (which project to the generator $\sigma$ of $\Z/n$) 
as 
$$\Theta_\pi(g\rtimes \sigma) = \Theta'({\rm Nm \,}g)\,\,\,\,\,\,\,\,\,\,\,\,\,\,\,\,\,\, ~~~~~~(*)$$
where $\Theta'$ is the restriction of the character of the representation 
$\pi_1 \otimes \pi_2 \otimes \cdots \otimes \pi_n$ of $\GL_m(\C)^n$ to the diagonally
embedded subgroup 
$$\GL_m(\C) \hookrightarrow \GL_m(\C)^n.$$ 

\vspace{4mm}
{\bf Remark 2:} The  kind of character relationship in $(*)$ is very familiar in  representation theory
of real and $p$-adic groups, and in fact was at the origin of this work 
which was to find aspects of the character theory of the generalized Steinberg representation $St_n(\pi)$ (or the generalized
Speh representation) for $\GL_{mn}(k)$, $k$ a $p$-adic field now,  which will allow one to retrieve $\pi$  from $\St_n(\pi)$. 
It  seems a 
worthwhile problem to find generalization of the theorem proved here for $\GL_{mn}(\C)$ to
$p$-adic $\GL_{mn}(k)$. Finding analogue of Kostant's theorem on character values at the 
Coxeter conjugacy class  seems an interesting question for $\GL_m(k)$ already. For instance, 
it appears that the character of a supercuspidal representation of $\GL_m(k)$ at the Coxeter
element is always zero, say in the tame case of $m$ coprime to the residue characteristic of $k$, since among other things, 
the Coxeter element cannot be transported to a division algebra.
The character relationship proved in the previous theorem reminds one of Jacquet modules but note that here we need only a Levi subgroup and its normalizer, rather 
than a parabolic, so offers more options for reductive groups over $p$-adic fields, such as 
$\U(m)^n\hookrightarrow \U(mn)$.

\vspace{ 4mm}

{\bf Example :} We calculate 
the character of irreducible representations ${\rm Sym }^{k}( \C^{2n}) $ and $\Lambda^k(\C^{2n})$ 
at the elements of $\GL_{2n}(\C)$ with eigenvalues $(t_1,\cdots, t_n, -t_1,\cdots, -t_n)$.

Since the character   of ${\rm Sym }^{k}( \C^{2n}) $ at the element of $\GL_{2n}(\C)$ with eigenvalues $(x_1,\cdots, x_{2n})$ is given by the
generating function,
$$\frac{1}{(1-tx_1)(1-tx_2) \cdots (1-tx_{2n})},$$
therefore the character of ${\rm Sym }^{k}( \C^{2n}) $ at the element of $\GL_{2n}(\C)$ with eigenvalues 
$(t_1,\cdots, t_n, -t_1,\cdots, -t_n)$ 
is given by the 
generating function,
$$\frac{1}{(1-tx_1) \cdots (1-tx_n) (1+tx_1) \cdots (1+tx_n)} = \frac{1}{(1-t^2x^2_1) \cdots (1-t^2x^2_n)}.$$

Thus the 
character of ${\rm Sym }^{k}( \C^{2n}) $ at the element of $\GL_{2n}(\C)$ with eigenvalues $(t_1,\cdots, t_n, -t_1,\cdots, -t_n)$ 
is nonzero only if $k$ is even, say $k= 2 \ell$, in which case the character at 
the element of $\GL_{2n}(\C)$ with eigenvalues $(t_1,\cdots, t_n, -t_1,\cdots, -t_n)$ is the 
character of ${\rm Sym}^{\ell}(\C^n)$ at the element of $\GL_n(\C)$ with eigenvalues $(t_1^2,\cdots, t_n^2)$.

Similarly, the character of $\Lambda^{k}( \C^{2n}) $ at the element of $\GL_{2n}(\C)$ with eigenvalues $(t_1,\cdots, t_n, -t_1,\cdots, -t_n)$ 
is nonzero only if $k$ is even, say $k= 2 \ell$, in which case the character at 
the element of $\GL_{2n}(\C)$ with eigenvalues $(t_1,\cdots, t_n, -t_1,\cdots, -t_n)$ is the 
character (up to a sign of $(-1)^\ell$) of $\Lambda^{\ell}(\C^n)$ at the element of $\GL_n(\C)$ with eigenvalues $(t_1^2,\cdots, t_n^2)$.

In both the examples, ${\rm Sym }^{k}( \C^{2n}) $ as well as $\Lambda^k(\C^{2n})$, the irreducible representation of the Levi subgroup
which is $\GL_n(\C) \times \GL_n(\C)$ is trivial on one of the factors, which is of course in accordance with Theorem 2. 
\section{A reformulation}

As we look for possible generalizations of Theorem 2, we rephrase this theorem  using  highest co-weights  instead of
highest weights which we recall now.

For any integer $a \geq 1$, 
let $\underline{\delta}_{a}$ be defined to be the $a$-tuple of integers,
$$a-1 \geq a-2 \geq \cdots \geq 1 \geq 0.$$

Interchanging character and co-character groups of a torus $T$ 
introduces a contravariant functor $T \rightarrow \wT$ such that,
$$\widehat{T}(\C) = {\rm Hom} [\C^\times, \widehat{T}] \otimes_\Z  \C^\times  = 
{\rm Hom}[T,\C^\times]\otimes_\Z \C^\times. $$ 

Thus given a 
character $\chi: T \rightarrow \C^\times$, it gives rise to a co-character,
 $$\widehat{\chi}: \C^\times \longrightarrow  \widehat{T}(\C)= {\rm Hom}[T,\C^\times]\otimes_\Z \C^\times, $$ 
 given by $z \longrightarrow \chi \otimes z$.

One can rephrase the theory of highest weight representations of $G(\C)$
to say that there is a natural bijective correspondence between finite dimensional 
irreducible representations of $G(\C)$ and conjugacy classes of algebraic homomorphisms:
$$\phi_{\lambda}: \C^\times \longrightarrow \widehat{G}(\C).$$

It is customary to absorb the factor $\rho$,  half the sum of positive 
roots, in $\lambda$ itself in this Harish-Chandra-Langlands parametrization, so $\phi_{\lambda} = \widehat{\lambda\cdot \rho}$.
This process eliminates some of the 
cumbersome $\rho$  shifts, in  particular  in the Weyl character formula. 
These characters $\lambda$ which index 
$\phi_\lambda$ are therefore 
strictly dominant, and $\phi_\lambda$ defines a natural Borel subgroup of $\widehat{G}$ 
(consisting of weights for $\C^\times$ which are non-negative for the adjoint action
of $\C^\times$  on the Lie algebra of $\widehat{G}(\C)$ 
via $\phi_{\lambda}$). 
We will fix the maximal torus $\widehat{T}$ of $\widehat{G}$ as the centralizer
of $\phi_\lambda$, and also fix the Borel subgroup $\widehat{B}$ of $\widehat{G}$ containing $\widehat{T}$
as just defined. 

We recall however that $\rho$ being half the sum of roots of $T$ inside $B$ is not necessarily a 
character of $T$. To make sense of $\rho$ we will have to restrict ourselves to $G$ a semi-simple simply
connected group in which case $\wrho$ is a cocharacter $\wrho: \C^\times \rightarrow \widehat{T} <  \widehat{G}$.   
In our case of $\GL_a(\C)$, where $\rho$ does not make sense, we redefine $\rho$ to be $\rho_a$ as above.

Here is the reformulation of Theorem 2 in this language.

\begin{thm}
Let $(t,c)$ be an element in $ \GL_m(\C) \times \GL_n(\C)$ with $c_n$ the Coxeter 
conjugacy class in $\GL_n(\C)$, and consider the image $t \cdot c_n$ of $(t,c_n)$ 
inside $\GL_{mn}(\C)$
under the natural map $\GL_m(\C) \times \GL_n(\C) \rightarrow \GL_{mn}(\C)$. Then the
character of a finite dimensional irreducible representation $\pi_\lambda$ of $\GL_{mn}(\C)$
corresponding to the cocharacter $\phi_\lambda: \C^\times \longrightarrow \GL_{mn}(\C)$, 
takes nonzero value at some element
of the form $t\cdot c_n$ if and only if  
$\phi_\lambda(e^{2\pi i /n})$
belongs to  the  conjugacy class in $\GL_{mn}(\C)$ defined by $1 \cdot c_n$.
Assume after 
conjugation that  $\phi_\lambda(e^{2\pi i /n})= 1 \cdot c_n$. 
If this is the case, then the image $\phi_\lambda( \C^\times)$  of $\phi_\lambda$ under 
$\phi_\lambda: \C^\times \longrightarrow \GL_{mn}(\C)$ commutes 
 with $\phi_\lambda(e^{2\pi i /n})= 1 \cdot c_n$, and 
therefore belongs to $\GL_m(\C)^n \hookrightarrow \GL_{mn}(\C)$ which is the commutant of $1\cdot c_n$ in $\GL_{mn}(\C)$.  
Let $\phi_{0,n} = \wrho_n: \C^\times \rightarrow \GL_n(\C)$ be the cocharacter corresponding to the trivial representation of $\GL_n(\C)$,
so that $\wrho_n(e^{2\pi i /n}) = c_n$. Thus,  $\phi_\lambda \wrho_n^{-1}: \C^\times \longrightarrow \GL_m(\C)^n \hookrightarrow \GL_{mn}(\C)$ 
is trivial on $e^{2 \pi i /n}$, and therefore there is a cocharacter $\phi_\mu$ with $\phi_\lambda \wrho_n^{-1} = \phi_\mu^n: \C^\times \longrightarrow \GL_m(\C)^n$.
 The coweight 
$\phi_\mu :\C^\times \longrightarrow \GL_m(\C)^n $  (which is regular inside $\GL_m(\C)^n$ although not necessarily in $\GL_{mn}(\C)$), defines a finite dimensional irreducible representation $\pi_\mu$ of $\GL_m(\C)^n$ such that the character of $\pi_\lambda$ at $t \cdot c_n$ is the 
same as the character of $\pi_\mu$ at the 
 element   $t^n$ inside $\Delta \GL_m(\C) \hookrightarrow \GL_m(\C)^n$. 
\end{thm}

\section{Generalization}

The author has considered several  possible options for generalizing the theorem of this paper to other groups. First let us  recapitulate the context
and what we would like to see happen. Suppose $G$ is a connected reductive group over $\C$ with $M$ a reductive subgroup such that for the connected 
component
$M_0$ of $M$, $M/M_0$ is a finite cyclic group. Write $M=M_0\cdot \langle \theta \rangle$. Note that we have a split exact sequence,
$$1 \rightarrow {\rm Int}(M_0)\rightarrow {\rm Aut}(M_0) \rightarrow {\rm Out}(M_0)\rightarrow 1,$$
where ${\rm Aut}(M_0)$ is the group of automorphisms of $M_0$, and ${\rm Int}(M_0)$ denotes the subgroup of  ${\rm Aut}(M_0)$ which consists of inner automorphisms
from elements of $M_0$, with splittings $: {\rm Out}(M_0) \longrightarrow {\rm Aut}(M_0)$ given by {\it pinnings} on $M_0$.  This allows us to  assume that
$M=M_0\cdot \langle \theta \rangle$ 
such that the automorphism induced by the action of  $\theta$ on $M_0$ preserves  a pinning on $M_0$. It seems 
un-necessary to assume that $M=M_0\cdot \langle \theta \rangle$ is actually a semi-direct product, as it is indeed not the case even for $M$ the
normalizer of the split torus in $\SL_2(\C)$.

Observe that the conjugacy class of the element $m\cdot \theta$ 
in $M$ is the same as the $\theta$-conjugacy class in $M_0$ (where $m$ in $M$ is conjugate to $n \cdot m \cdot \theta(n^{-1})$).

Assume that the 
$\theta$-conjugacy class in $M_0$ are the same as conjugacy class in another group $M'$, and that there is a natural map from 
conjugacy classes in $M'$ to conjugacy classes in $M_0$. Thus there is a map from $\theta$-conjugacy classes in $M_0$ to conjugacy classes in $M_0$. Call this map 
a {\it Norm} map ${\rm Nm}$ from $\theta$-conjugacy classes in $M_0$ to conjugacy classes in $M_0$.
The aim then is to express the character of an irreducible
representation of $G$ at an element of $M$ of the form $m\cdot \theta$ as the character of an  irreducible representation of $M$ at ${\rm Nm}(m)$.
In the case of $\GL_{mn}(\C)$, from an irreducible representation of $\GL_{mn}(\C)$ which is given by an 
unordered $mn$-tuple of characters of $\C^\times$, we need to construct an irreducible representation of $M_0= \GL_m(\C)^n$ 
(well defined only up to permutation of the factors in $\GL_m(\C)^n$, and multiplication of 
one dimensional characters on each factor with product trivial) which   amounts to having
unordered $n$-tuples of unordered $m$-tuples. 
This construction which we would like to emphasize is not 
obvious  needs to be done in some generality, and arises it seems
by looking at centralizers of certain elements in the dual group of $G$, as in Theorem 3. 

For character relations as desired in the previous paragraph to happen, it seems necessary that,

\begin{enumerate}
\item $M$ is not contained in a proper connected subgroup of $G$.

\item $M_0$ has no non-trivial normal semi-simple subgroup on which $\theta$ operates trivially. (We have already assumed that $\theta$ preserves a 
pinning of $M_0$.)
\end{enumerate}

Given the theorem for $\GL_n(\C)$, the most obvious next case to consider will be Levi subgroups
of parabolics in reductive groups for which there are no other associate parabolics, so that the normalizer of the Levi subgroup in the ambient 
reductive group divided by the Levi subgroup is a Weyl group acting irreducibly on the center of the Levi subgroup modulo center of the group.

The other option is to look for centralizer of powers of the Coxeter element in $G$. 

Neither of the two options seems to be working out (as they do not satisfy condition 2 above in general), and the best hope seems the following very specific situation.

Let $G= \Sp_{2mn}(\C)$ containing $M_0=\GL_m(\C)^n$ as a Levi subgroup, sitting as the block diagonal matrices,

$$\left  ( \begin{array}{cccccccc}
g_1& {} & {} & {}  &  {}   & {}  &{} &  \\
{}& g_2  & {} & {}   &  {}  &  {}  & {} &   \\
{}& {}  & \cdot  & {}  &   {} & {}   & {}   &\\
{}& {}  &   & {g_n}  &   {} & {}   &{} &  \\
{}& {}  &   & {}  &   {g_1^{*}} & {}   &{}   &\\
{} &  {}  &   {}   &  &  & g_2^{*} &    & \\
{}& {}  &   & {}  &   {} &  & {\cdot}      &    \\
{} & {}   &  {}  & {}  &  & & & g_n^{*}  
\end{array} \right)\cdot
$$
with $g^* = j ^tg^{-1} j^{-1}$ where $j$ is the anti-diagonal matrix with alternating 1 and $-1$. Define $\theta$ to be the cyclic permutation:

\begin{eqnarray*}
\left  ( \begin{array}{cccccccc}
g_1& {} & {} & {}  &  {}   & {}  &{} &  \\
{}& g_2  & {} & {}   &  {}  &  {}  & {} &   \\
{}& {}  & \cdot  & {}  &   {} & {}   & {}   &\\
{}& {}  &   & {g_n}  &   {} & {}   &{} &  \\
{}& {}  &   & {}  &   {g_1^{*}} & {}   &{}   &\\
{} &  {}  &   {}   &  &  & g_2^{*} &    & \\
{}& {}  &   & {}  &   {} &   & {\cdot}      &    \\
{} & {}   &  {}  & {}  &  & & & g_n^{*} 
\end{array} \right)
& \longrightarrow & 
\left  ( \begin{array}{cccccccc}
g_n^{*}& {} & {} & {}  &  {}   & {}  &{} &  \\
{}& g_1  & {} & {}   &  {}  &  {}  & {} &   \\
{}& {}  & \cdot  & {}  &   {} & {}   & {}   &\\
{}& {}  &   & {g_{n-1}}  &   {} & {}   &{} &  \\
{}& {}  &   & {}  &   {g_n} & {}   &{}   &\\
{} &  {}  &   {}   &  &  & g_1^{*} &    & \\
{}& {}  &   & {}  &   {} &   & {\cdot}      &    \\
{} & {}   &  {}  & {}  &  & & & g_{n-1}^{*} 
\end{array} \right),
\end{eqnarray*}
which corresponds to 
$$\theta: (g_1,g_2,\cdots,g_{n-1}, g_n) \longrightarrow (g_n^{*},g_1,g_2,\cdots, g_{n-1}),$$
and $M= M_0 \cdot \langle \theta \rangle$,
but even this example we have not worked out. There are of course very similar examples to consider for $\SO_{2mn+1}(\C)$ as well as $O_{2mn}(\C)$.

\vspace{4mm}

{\bf Example :} One of the simplest contexts for generalization  would be $\GL_n(\C) \hookrightarrow \Sp_{2n}(\C)$, sitting as 
a Levi subgroup of a Siegel parabolic. In this case, $\GL_n(\C)$ has a nontrivial normalizer $N\GL_n(\C)$
inside $\Sp_{2n}(\C)$ with the exact sequence,
$$1 \rightarrow \GL_n(\C) \rightarrow N\GL_n(\C) \rightarrow \Z/2 \rightarrow 1.$$

We calculate 
the character of ${\rm Sym }^{k}( \C^{2n}) $  --- which is known to be an irreducible representation 
of $\Sp_{2n}(\C)$ 
--- at the elements of $\Sp_{2n}(\C)$ which belong to $N\GL_n(\C) = \GL_n(\C) \cdot \theta$ but not to $\GL_n(\C)$ where $\theta$ is an element of $\Sp_{2n}(\C)$ 
which operates on $\GL_n(\C)$ by the outer automorphism $g \rightarrow J {}^tg^{-1}J^{-1}$ where $J$ is the anti-diagonal matrix with entries
which are alternatingly 1 and $-1$.

 It can be seen that for $n=2m$, any element of $\GL_n(\C)\cdot \theta$ can be represented as $t \cdot \theta$ where $t=(t_1,\cdots, t_m, 
t_m^{-1},\cdots, t_1^{-1}) \cdot \theta$, and such elements are conjugate inside $\GL_{4m}(\C)$ to 
$(t_1,\cdots, t_m, t_1^{-1}, \cdots t_m^{-1}, -t_1,\cdots, -t_m, -t_1^{-1},\cdots, -t_m^{-1})$.

Since the character   of ${\rm Sym }^{k}( \C^{2n}) $ at the element of $\GL_{2n}(\C)$ with eigenvalues $(x_1,\cdots, x_{2n})$ is given by the
generating function,
$$\frac{1}{(1-tx_1)(1-tx_2) \cdots (1-tx_{2n})},$$
therefore the character of ${\rm Sym }^{k}( \C^{4m}) $ at the element of $\GL_{4m}(\C)$ with eigenvalues 
$(t_1,\cdots, t_m, t_1^{-1}, \cdots t_m^{-1}, -t_1,\cdots, -t_m, -t_1^{-1},\cdots, -t_m^{-1})$ 
is given by the 
generating function,
$$\frac{1}{(1-t^2t^2_1) \cdots (1-t^2t^2_m)(1-t^2t^{-2}_1) \cdots (1-t^2t^{-2}_m)} 
$$

Thus the 
character of ${\rm Sym }^{k}( \C^{4m}) $ at the element of $\GL_{4m}(\C)$ with eigenvalues 
$(t_1,\cdots, t_m, t_1^{-1}, \cdots t_m^{-1}, -t_1,\cdots, -t_m, -t_1^{-1},\cdots, -t_m^{-1})$  
is nonzero only if $k$ is even, say $k= 2 \ell$, in which case the character at 
the element of $\GL_{4m}(\C)$ with eigenvalues 
$(t_1,\cdots, t_m, t_1^{-1}, \cdots t_m^{-1}, -t_1,\cdots, -t_m, -t_1^{-1},\cdots, -t_m^{-1})$ is the 
character of ${\rm Sym}^{\ell}(\C^{2m})$ at the element of $\GL_{2m}(\C)$ with eigenvalues $(t_1^2,\cdots, t_m^2, t_1^{-2},\cdots, t_m^{-2})$.

\vspace{6mm}

{\bf Acknowledgement: } The present work was conceived during the summer of 2008 while the author was
enjoying the hospitality of the University of California at  San Diego hosted by Wee Teck Gan, and partially supported 
by the Clay Math Institute.  The author thanks Prof.\!\! Wee Teck Gan for his hospitality, and UCSD as well as 
CMI for their  support.

\vspace{1cm}

{\bf \Large Bibliography}

\vspace{1cm}

[Ko] Kostant, Bertram:
{\it On Macdonald's $\eta$-function formula, the Laplacian and generalized exponents.}
Advances in Math. 20 (1976), no. 2, 179--212. 

\vspace{3mm}

[KLP] Kumar, Shrawan, Lusztig, George and Prasad, Dipendra: {\it Characters of simplylaced nonconnected
groups versus characters of nonsimplylaced connected groups}. Contemporary 
Math., vol 478, AMS, pp. 99-101.

\vspace{3mm}

[La] Langlands, Robert:{\it  Base change for GL(2).} 
Annals of Mathematics Studies, 96. Princeton University Press, Princeton, N.J.; 
University of Tokyo Press, Tokyo, 1980.

\end{document}